\documentclass[12pt, reqno, letterpaper, titlepage, final]{article}

\usepackage{amsmath}
\usepackage{amssymb}
\usepackage{amsthm}
\usepackage[english]{babel}
\usepackage{color}
\usepackage{delarray}
\usepackage{enumerate}
\usepackage{epsfig}
\usepackage{graphicx}
\usepackage{latexsym}
\usepackage{times}
\usepackage{cite}

\theoremstyle{plain}
\newtheorem{theorem}{Theorem}[section]

\newtheorem{lemma}{Lemma}[section]
\newtheorem{proposition}{Proposition}[section]
\theoremstyle{definition}

\newtheorem{remark}{Remark}[section]

\makeatletter
\renewcommand\@biblabel[1]{#1.}
\makeatother

\newcommand{\R}{\mathbb{R}}

\newcommand{\norm}[1]{\|#1\|}

\begin{document}

\centerline{\Large{\bf Decomposition algorithms }}
\centerline{\Large{\bf for globally solving mathematical programs}}
\centerline{\Large{\bf with affine equilibrium constraints \footnote {This work is supported by the Vietnam National Foundation for Science Technology Development (NAFOSTED).}}}

\vskip 0.8cm
\centerline{ {\bf L. D. Muu}\footnote{Institute of Mathematics, VAST, Hanoi, Vietnam, e-mail: ldmuu@math.ac.vn}~$\bullet$~ {\bf T. D. Quoc} \footnote{Hanoi University of Science, Hanoi, Vietnam.} ~$\bullet$~ {\bf L. T. H. An}\footnote{Metz University, France.}~$\bullet$~ {\bf P. D. Tao} \footnote{INSA, Rouen, France.}}

\vskip 1.0cm

{\leftskip 1.2cm {\noindent \bf{Abstract }}{ 
A mathematical programming problem with affine equilibrium constraints (AMPEC) is a bilevel programming problem where the lower one is a parametric affine variational inequality. We formulate some classes of bilevel programming in forms of AMPEC. 
Then we use a regularization technique to formulate the resulting problem as a mathematical program with an additional constraint defined by the difference of two convex functions (DC function). A main feature of this DC decomposition is that the second component depends upon only the parameter in the lower problem. This property allows us to develop branch-and-bound algorithms for globally solving AMPEC where the adaptive rectangular bisection takes place only in the space of the parameter.
As an example, we use the proposed algorithm to solve a bilevel Nash-Cournot equilibrium market model. Computational results show the efficiency of the proposed algorithm.
}}


\bigskip


\section{Introduction}\label{sec:intro}
We consider the following mathematical programming problem with affine (not necessarily monotone) variational inequality constraints, that we call shortly AMPEC:
\begin{align}\label{eq:mpec_prob}
\min_{x\in\R^n,y\in\R^m} &f(x,y) \tag{P}\\
\textrm{s.t.}~~~~~&(x, y) \in S, \label{eq:mpec_c1}\\
&x\in C, ~ (Ax+By+a)^T(v-x) \geq 0, ~\forall v\in C, \label{eq:mpec_c2} 
\end{align}
where $\emptyset \neq S \subseteq \R^{n+m}$, $\emptyset \neq C \subseteq \R^n$ are two closed convex sets, $f :\R^{m+n} \to \R$ is a convex function, $A$, $B$ are given appropriate real matrices, and $a\in\mathbf{R}^n$.
This class of optimization problems is known to be very difficult to solve due to its nonconvexity, nondifferentiability and loss of constraint qualification. However such problems arise frequently in applications, for example, in shape optimization, design transportation network, economic modeling and data mining. 
A natural way to handle a nested problem such as Problem \eqref{eq:mpec_prob} is to reduce it into an one-level optimization problem by using the Krush-Kuhn-Tucker theorem for the lower variational inequality.
Several algorithms for globally solving the reduced mathematical programs with complementarity constraints are proposed (see, e.g.
\cite{An2009, Anh2006, Muu2007, Thoai2005}).
Since the number of the complementarity constraints is just equal to the number of constraints defining the set $C$ in the lower variational inequality problem, these global optimization algorithms become expensive when the number of constraints is high, for example, when $C := \{ x \in R^n~|~x\geq 0,~ c_j(x) \leq 0, ~j = 1,\dots,p \}$ (often appears in practice) with either $n$ or $p$ are somewhat large.
 
In this paper, we propose another solution-approach to AMPEC without using the Krush-Kuhn-Tucker theorem for the lower variational inequality. Instead, we use a regularization technique to formulate AMEC as a mathematical program with an additional constraint defined by $g_1(x,y) - h_1(x,y) \leq 0$, where $g_1$ and $h_1$ are differentiable convex functions. The main feature of this constraint is that the second component $h_1$ can be chosen such a way so that it depends upon only the parameter $y$. Moreover, in some special important cases such as bilevel convex quadratic problems, $h_1$ is separable. This formulation allows us to develop decomposition branch-and-bound algorithms for globally solving AMPEC where the branching operation involving only the parameter in the lower variational inequality. Unlike the existing global optimization algorithms mentioned above, the proposed algorithms can solve AMPEC where the constraint set $C$ is given as $C := \{ x \in R^n~|~x\geq 0, ~c_j(x) \leq 0, ~j =1,\dots,p\}$ with $n$ and $p$ relatively large.
As an example, we use the proposed algorithm to find a global optimal equilibrium pair to a bilevel Nash-Cournot equilibrium market model.
We tested the proposed algorithm by some randomly generated data. The numerical results show that our algorithm can solve this bilevel model with high dimension.
 
The paper is organized as follows. In the next section we give DC formulations to AMPEC by using suitable regularization matrices. Some important special cases of AMPEC are presented at the end of this section. The third section is devoted to description of a branch-and-bound algorithm for globally solving a bilevel Nash-Cournot equilibrium market model by using a DC decomposition, where the second component is separable and depends upon only the parameter $y$. We close the paper with some computational experiences and results. 

\section {DC Formulations and Examples}\label{sec:DC_formulation}
In Problem \eqref{eq:mpec_prob}, as usual, we will refer to $x$ as a primary variable or decision variable and $y$ as a parameter. We call $(x,y)$ a feasible point to \eqref{eq:mpec_prob} if $(x,y) \in S$ and $x$ solves the lower  variational inequality \eqref{eq:mpec_c2}. Note that when $A$ is  symmetric positive semidefinite, the variational inequality \eqref{eq:mpec_c2} is equivalent to the parametric convex quadratic problem
\begin{equation}\label{eq:quadratic_prob}
\min \big\{ \varphi (x,y):=\frac{1}{2} x^TAx + (By + a)^T x ~|~ x \in C \big\}.
\end{equation}
In this case, Problem \eqref{eq:mpec_prob} becomes a bilevel convex program
\begin{equation}\label{eq:bilevel_prob}
\min \{f(x,y) ~|~ (x,y) \in S\} \tag{BP},
\end{equation}
where $x$ solves the convex quadratic program
\begin{equation}\label{eq:bilevel_quad_prob}
\min\big\{\varphi (x,y) := \frac{1}{2}x^TAx + (By + a)^T x ~|~ x\in C \big\}.
\end{equation}
In the general case, when $A$ is indefinite, the variational inequality \eqref{eq:mpec_c2} is not necessarily equivalent to the problem \eqref{eq:bilevel_quad_prob}. So Problem \eqref{eq:mpec_prob}, in general, cannot be reformulated as a bilevel problem of the form \eqref{eq:bilevel_prob} \cite{Dempe2003, Lee2005}.

\subsection{DC Formulations}\label{subsec:DF_form}
The main difficulty of Problem \eqref{eq:mpec_prob} is that the constraint defined by the variational inequality \eqref{eq:mpec_c2} is neither convex nor given explicitly as a constrained set of an ordinary mathematical programming problem. A natural way is to reduce Problem \eqref{eq:mpec_prob} to an ordinary mathematical programming problem. We will reformulate Problem \eqref{eq:mpec_prob} as a smoothly DC program. For this purpose, we use a gap function introduced in \cite{Taij1996} to formulate the variational inequality \eqref{eq:mpec_c2} as an equation defined by a smoothly DC function. We recall that a function $f$ is said to be DC on a convex set $D$ if it can be expressed as the difference of two convex functions on $D$, i.e. $f=g-h$, where $g$ and $h$ are convex on $D$.

More precisely, for each $(x,y)$, following the idea from \cite{Taij1996} we define the function $g(x,y)$ by setting
\begin{equation}\label{eq:gap_func}
g(x,y) := \max_{v\in C} \big\{(x-v)^T(Ax + By + a) -\frac{1}{2}(v-x)^T G(v-x) \big\},
\end{equation}
where $G$ is an arbitrary symmetric positive definite $(n\times n)$-matrix. We refer to $G$ as a {\it regularization matrix}. Since $G$ is positive definite, the problem defining $g(x,y)$ is uniquely solvable for every $(x,y)$, i.e. $g$ is well-defined.

The following lemma gives the properties of the gap function $g$ whose proof can be done similarly to the proof of Lemma 2.1 in \cite{Taij1996}.

\begin{lemma}\label{le:gap_func}\ Let $g$ be given by \eqref{eq:gap_func}. Then
\begin{itemize}
\item[$\textrm{(i)}$] $g(x,y) \geq 0$ for every $(x,y) \in C\times \R^m$,
\item[$\textrm{(ii)}$] $(x,y) \in S$, $x\in C$, $g(x,y) = 0$ if and only if $(x,y)$ is feasible solution to Problem \eqref{eq:mpec_prob}.
\end{itemize}
\end{lemma}
The following proposition shows that, with a suitable choice of the regularization matrix $G$, the function $g$ can be decomposed as the difference of two convex functions (DC function). Note that any symmetric matrix $A$ can be expressed as $A = A_1- A_2$, where $A_1$ is symmetric positive definite and $A_2$ is symmetric. In what follows by $\text{diag}(\alpha) $ we denote the diagonal matrix whose every diagonal entry is $\alpha$.

\begin{proposition}\label{pro:dc_decomp1}
Suppose that $A$ is symmetric and $A = A_1 - A_2$, where $A_1$ is a symmetric positive definite matrix and $A_2$ is a symmetric matrix such that $A_2 + \frac{1}{2}U^TU$ is positive (semi)definite, and $U$, $V$ are two appropriate matrices satisfying $U^TV=B$.
Let $G = 2A_1$. Then
\begin{equation}\label{eq:dc_decomp1}
g(x,y) = g_1(x,y) - h_1(x,y), 
\end{equation}
where $g_1$ and $h_1$ are two differentiable convex functions given by
\begin{align}\label{eq:g1_func}
g_1(x,y) & = \frac{1}{2}\norm{Ux + Vy}^2 +a^Tx \notag\\
&+ \max_{v\in C}\big\{ [(A_1+A_2)x -By - a]^Tv - v^TA_1v \big\},
\end{align}
and
\begin{equation}\label{eq:h1_func}
h_1(x,y) = \frac{1}{2}x^T(2A_2+U^TU)x + \frac{1}{2}\norm{Vy}^2.
\end{equation} 
\end{proposition}

\begin{proof}
With a simple arrangement from \eqref{eq:gap_func}, it shows that
\begin{align}\label{eq:g_func_exp}
g( x, y ) &= x^TAx - \frac {1}{2}x^TGx + x^TBy + a^Tx \notag\\
&+ \max_{v\in C} \{ -v^TAx - v^TBy - a^Tv - \frac{1}{2}v^TGv + x^TGv\}. 
\end{align}
Since $A = A_1 - A_2$ and $G = 2A_1$, the last expression implies
\begin{align}\label{eq:gap_express}
g(x, y) &= -x^TA_2x + x^TBy + a^Tx \notag\\
& + \max_{v\in C}\{ - v^TA_1v + [(A_1+A_2)x-By-a]^Tv \}.
\end{align}
On the other hand, since $B=U^TV$ we can express
\begin{equation*}\label{eq:xBx_exp}
2 x^T B y = 2x^TU^TVy = \norm{ Ux + Vy }^2 - \norm{Ux}^2 - \norm{Vy}^2 .
\end{equation*}
Substituting this expression into \eqref{eq:gap_express} we get
\begin{align*}
g( x, y) &= \frac{1}{2}\norm{Ux + Vy}^2 - \frac{1}{2}\norm{Ux}^2 - \frac{1}{2}\norm{Vy}^2 +a^Tx\\
&+ \max_{v\in C}\{ - v^TA_1v + [(A_1+A_2)x - By - a]^Tv \}. 
\end{align*}
Hence,
\begin{equation*}
g(x,y) = g_1(x,y) - h_1(x,y), 
\end{equation*}
where $g_1$ and $h_1$ are two functions given by \eqref{eq:g1_func} and \eqref{eq:h1_func}, respectively. Since $A_2 + \frac{1}{2}U^TU$ is positive semidefinite, $h_1$ is convex. Clearly, $h_1$ is differentiable everywhere, while $g_1$ is differentiable everywhere because the convex program (strongly quadratic concave maximization):
\begin{equation*}
\max_{v\in C}\big\{ -v^TA_1v + [(A_1+A_2)x - By - a]^Tv\big\} 
\end{equation*}
is uniquely solvable for any $(x,y)$. 
\end{proof}

\begin{remark}\label{re:derivative_of_g} 
From (\ref{eq:g1_func}), by a simple computation we have
\begin{align}
&\nabla_x g_1(x, y ) = U^T(Ux + Vy) + a + (A_1+A_2)^Tz(x, y), \label{eq:g1_derivative_x}\\
&\nabla_y g_1(x. y) = V^T(Ux + Vy) - B^Tz(x, y), \label{eq:g1_derivative_y}
\end{align}
where $z(x,y)$ is a unique solution of the strongly convex quadratic program
\begin{equation*}
\max_{v\in C}\big\{ -v^TA_1v + [(A_1+A_2)x - By - a]^Tv \big\}.
\end{equation*}
\end{remark}
 
\begin{remark}\label{re:B_expression}
Since matrices $U$ and $V$ in Proposition \ref{pro:dc_decomp1} can be arbitrary, we can choose $U$ and $V$ such that $V$ has a simple form. 
For example, if we choose $U = [(\Sigma B^{+})^T]^{+}$, where $B^{+}$ is the (Moore-Penrose) pseudo-inverse of $B$ and $\Sigma$ is a diagonal matrix, then $V$ is a diagonal matrix, precisely, $V=\Sigma$.
\end{remark}
We call the DC decomposition $g(x,y) = g_1(x,y) - h_1(x,y)$, where $g_1$ and $h_1$ are given by \eqref{eq:g1_func} and \eqref{eq:h1_func} respectively, a {\it spectral decomposition}. In this decomposition, the function $h_1$ is a quadratic form, even separable quadratic if $2A_2+U^TU$ is diagonal. The separable quadratic property of $h_1$ is useful when applying to global algorithms that use the convex envelope of $-h_1$ (see Section \ref{sec:global_algorithm} below). 
 
Using Proposition \ref{pro:dc_decomp1}, Problem \eqref{eq:mpec_prob} is reformulated equivalently to a DC constrained optimization problem of the form
\begin{align}\label{eq:mpec2_prob}
&\min_{x\in\R^n,y\in\R^m} f(x,y) \tag{$\textrm{P}_1$}\\
&\textrm{s.t.} ~~ (x,y) \in S,~ x\in C \label{eq:convex_con}\\
& g(x,y) = g_1(x,y) - h_1(x,y) \leq 0, \label{eq:dc_con} 
\end{align}
where $g_1$ and $h_1$ are given by \eqref{eq:g1_func} and \eqref{eq:h1_func}, respectively.
 
Formulation \eqref{eq:mpec2_prob} allows that theory and methods in smooth and DC optimization both global and local can be applied to mathematical programs with affine equilibrium constraints.

\subsection{Special Cases}\label{subsec:special_cases}
In this subsection, we consider some special, but important, cases of Problem \eqref{eq:mpec_prob} and their reformulation in the form of  \eqref{eq:mpec2_prob}.

\subsubsection{Linear program with linear complementarity constraints}
\label{subsubsec:LPwithCC}
Note that when $C = \R^n_+$, $S$ is a polyhedron defined by
\begin{equation*}
S := \big\{(x,y) ~:~ Ax + By + a \geq 0 \big\}, 
\end{equation*}
and $f(x,y) = c^Tx + c^Ty$, Problem \eqref{eq:mpec_prob} becomes a linear program with an additional linear complementarity constraint of the form 
\begin{align}
&\min_{(x,y)} f(x, y), \label{eq:comp_prob}\tag{CP} \\
&\textrm{s.t.} ~~~x \geq 0,~ Ax + By + a \geq 0,~ x^T(Ax + By + a) = 0.\label{eq:CC}
\end{align}
For this program, the following gap function has been used \cite{An1999, Mangasarian1997, An2009}:
\begin{equation*}
p(x,y) = \sum_{j=1}^n\min\{x_j, (Ax + By + a)_j\} 
\end{equation*}
It has been shown that if $f$ is bounded from below, then there exists $t_{*} > 0$ such that for every $t \geq t_{*}$, Problem \eqref{eq:comp_prob} is equivalent to the following concave minimization problem
\begin{align*}
&\min_{(x,y)} \big\{ f(x,y) + tp(x,y)\big\}\\
&\textrm{s.t.}~~ x\geq 0,~ Ax + By + a \geq 0, 
\end{align*}
in the sense that their solution-set coincide. In \cite{Mangasarian1997} Mangasarian and Pang replaced $p$ by the differentiable function
\begin{equation*}
\min\{\sum_{j=1}^n r_jx_j + s_j(Ax + By + a)_j ~|~ r_j, s_j \geq 0, r_j + s_j = 1, j=1,\dots,n\}. 
\end{equation*} 
Note that the DC function $g(x,y) = g_1(x,y) - h_1(x,y)$ with $g_1$ and $h_1$ given as in Proposition \ref{pro:dc_decomp1} is a differentiable merit DC functions for \eqref{eq:comp_prob} without introducing $2n$-extra variables $r$ and $s$.
 
\subsubsection{Linear optimization over the Pareto-efficient set}
\label{subsubsec:LPoverPareto} 
Let $X \subset \R^n$ be a nonempty bounded polyhedron and $W$ be a $(p\times n)$-real matrix. Consider the vector optimization problem of the form
\begin{equation}\label{eq:VP}
\min\{Wx ~|~ x\in X\}. 
\end{equation}
We recall that a point $x^{*}\in X$ is said to be an {\it efficient solution} or a {\it Pareto solution} to \eqref{eq:VP}, if whenever $x\in X, Wx \leq Wx^{*}$, then $Wx = Wx^{*}$. Let $E(W,X)$ denote the set of all efficient solutions to \eqref{eq:VP}. Consider the optimization over the efficient set
\begin{equation*}
\min \{f(x) ~|~ x\in E(W,X)\} \tag{PP},
\end{equation*}
where $f$ is real valued convex function on $\R^n$. This problem has some applications in decision making and recently has been studied in many research articles (see, e.g.\cite{An1996, An2003, Benson1984, Luo1996, Muu2000, Philip1972} and references therein).
Note that since the efficient set is rarely convex, this problem is a nonconvex optimization problem.

It has been shown in \cite {Philip1972} that one can find a simplex $Y$ in $\R^p$ such that a point $x^{*}$ is efficient for \eqref{eq:VP} if and only if there exists $y^{*}\in Y$ such that
\begin{equation*}
(W^Ty^{*})^T(x - x^{*}) \geq 0, ~~\forall x\in X. 
\end{equation*}
Thus the above optimization problem over the efficient set can be formulated as the mathematical program with affine equilibrium constraint of the form 
\begin{equation}\label{eq:VP_EP}
\min\big\{f(x) ~|~ (x,y) \in X\times Y,~ (W^Ty)^T(v - x) \geq 0, ~\forall v \in X \big\}. \tag{EP} 
\end{equation}
By this way, a point  $x^*$ is a global optimization to (PP) if  and only there exists $y^* \in Y$ such that 
$(x^*,y^*)$ is a global optimal solution to (EP).
The latter problem is of the form \eqref{eq:mpec_prob} with $S = X\times Y$, $C = X$ and $A = 0$, $B = W^T$, $a = 0$. 
Since $A = 0$, we can apply Proposition \ref{pro:dc_decomp1}, for example, with
$A_1 = A_2 = I$, where $I$ is the identity matrix. Since $B = W^T$, $a = 0$ 
from Proposition \ref{pro:dc_decomp1} we have $g(x) = g_1(x) - h_1(x)$ with
\begin{align*}
&g_1(x,y) =\frac{1}{2}\norm{Ux + Vy}^2 + a^Tx + \max_{v\in C}\big\{(2Ix -W^Ty)v - v^TA_1v\big\},\\
&h_1(x,y)= \frac{1}{2}x^T (2A_2 + U^TU)x +\frac{1}{2}\norm{Vy}^2,
\end{align*}
where $U^TV = W^T$.
Thus, by Lemma \ref{le:gap_func}, we can  formulate  (PP) as the following optimization problem with a DC constraint
\begin{equation*}
\min \{f(x) ~|~ (x,y) \in X \times Y,~ g_1(x,y) - h_1(x,y) \leq 0 \}. 
\end{equation*}

\subsubsection{A bilevel Nash-Cournot oligopolistic equilibrium market model}
\label{subsubsec:Nash_EP} 
Suppose that there are $n$-firms (sectors) that supply a homogeneous product whose price $p$ at each sector $j$ $(j=1,\dots,n)$ depends on total producing quantity and is given by
\begin{equation*}
p(\sum_{j=1}^n x_j)=\alpha -\beta\sum_{j=1}^n x_j, 
\end{equation*}
where $\alpha > 0$, $\beta > 0$ are given constants, $x_j$ is the quantity of goods supplied by firm $j$ that we have to determine.
Suppose further that, to produce the goods, the firms need $m$-different materials represented by a vector $y \in \R^m$. Let $y_i$ ($i=1,\dots,m$)   be the quantity of material $i$ needed to produce a unique of goods. Let $c_{ji}$ denote the price of a unit
material $i$ for firm $j$ ($i=1,\dots,m$, $j= 1,\dots,n)$. When $c_{ji} \leq 0$, it means that firm $j$ is encouraged to use material $i$; for example, it is a waste material. Assume that the cost of firm $j$ is given by
\begin{equation*}
h_j(x_j,y):=x_j\sum_{i=1}^mc_{ji} y_i+\delta_j,\ \ j=1,\dots,n, 
\end{equation*}
where $\delta_j \geq 0$ is the fixed charge cost at firm $j$. Then the utility function of firm $j$ can be given by
\begin{equation*}
u_j(x,y) := p(\sum_{i=1}^n x_i)x_j - h_j(x_j,y). 
\end{equation*}
Let
\begin{align*}
&Y_i:= \{y_i\ :\ 0 \leq y_i \leq \xi_i \} ~~ (i=1,\dots,m),\\
&X_j:= \{\tau\ :\ 0 \leq \tau \leq \eta_j \} ~~ (j=1,\dots,n), 
\end{align*}
where $\xi_i$ is the upper bound of the material $i$, and $\eta_j$ is the upper bound of the quantity of the goods can be produced by firm $j$.

Let
\begin{equation*}
Y:=Y_1 \dots\times Y_m, \ \ X= X_1\times\dots\times X_n
\end{equation*}
be the feasible (strategy)-sets of the model.

Given $y\in Y$, each firm $j$ seeks to find its producing quantity $x_j$ such that its benefit $u_j(x,y)$ is maximal.
However, a maximal policy for all firms altogether, in general, does not exist. So they agree with an equilibrium point in the
sense of Nash.

By definition, a vector $(x^{*}_1,\dots,x^{*}_n) \in X_1\times \dots\times X_n$ is said to be a (Nash) equilibrium point with
respect to $y^{*}\in Y$ if, for all $x_j \in X_j$ and $j$,
\begin{equation*}
u_j(x^{*}_1,\dots,x^{*}_{j-1},x_j,x^{*}_{j+1},\dots,x^{*}_n,y^{*}) \leq u_j(x^{*}_1,\dots,x^{*}_{j-1},x^{*}_j,x^{*}_{j+1},\dots,x^{*}_n,y^{*}).
\end{equation*}
We will refer to such a pair $(x^{*},y^{*})$ as an {\it equilibrium pair} of the model.

Besides the utility function of each firm, there is another cost function (leader's objective function) $f(x,y)$ depending on $y$ and the quantity $x$ of the goods. The problem needs to be solved is of finding an equilibrium pair that minimizes
leader's objective function over the set of all equilibrium pairs. We call such a pair $(x^{*}, y^{*})$ a {\it global  optimal  equilibrium pair} to the model.
This problem can be reformulated as a mathematical program with affine equilibrium constraints. To this end, let 
\begin{align*}
&H_j(x_j,y):= \nabla_{x_j} h_j(x_j,y)~~ (j=1,\dots,n), \\
&e := (1,\dots,1)^T,\\
&\sigma_x :=\sum_{j=1}^n x_j. 
\end{align*}
Applying Proposition 3.2.6 in \cite{Konnov2000} we see that a point $(x_1,\dots,x_n)$ is equilibrium with respect to $y$ if and only if
it is a solution to the variational inequality problem
\begin{align*}
\text{Find} ~~& x\in X~~\text{such that:}~~ F(x,y)^T(z-x) \geq 0, ~~\text{for all}~~ z\in X, 
\end{align*}
where $F(x,y)$ is $n$-dimensional vector function whose $j$-th component is defined by
\begin{equation}\label{eq:Fij_func}
F_j(x,y):=H_j(x,y)-p(\sigma_x)-\nabla p(\sigma_x)x_j.
\end{equation}
Using \eqref{eq:Fij_func} and the definition of $H_j(x,y)$ we have
\begin{equation*}
F_j(x,y) = \sum_{i=1}^m c_{ji}y_i -\alpha +\beta\sum_{k=1}^n x_k+ \beta x_j ~~ (j=1,\dots,n) 
\end{equation*}
Thus
\begin{equation*}
F(x,y) = Ax + By+a, 
\end{equation*}
where
\begin{equation}\label{eq:matrix_A}
A = \begin{bmatrix}
2\beta&\beta &\beta &\cdots&\beta\\
\beta &2\beta&\beta &\cdots&\beta\\
\cdots&\cdots&\cdots&\cdots&\cdots\\
\beta &\beta &\beta &\cdots&2\beta
\end{bmatrix}
\end{equation}
and $B$ is an $(n\times m)$ matrix (independent of $x$) whose $B_{ij}$ entry is 
\begin{equation}\label{eq:matrix_B}
B_{ji}= c_{ji},\ j=1,\dots,n,\ i=1,\dots,m,
\end{equation}
and
\begin{equation}\label{eq:vector_a}
a=(-\alpha,\dots,-\alpha)^T \in \R^n.
\end{equation}
Thus the problem needs to be solved takes the form
\begin{align*}
\min_{x,y}~~&f(x, y)\\
\textrm{s.t.}~~~ &y\in Y := Y_1\times\cdots\times Y_m, ~ x\in X := X_1\times\cdots\times X_n\\
&\textrm{where} ~x~ \textrm{solves the parametric variational inequality}\\
& (Ax+By+a)^T(v-x) \geq 0, ~~\forall v\in X, 
\end{align*}
with $A$, $B$ and $a$ being given by (\eqref{eq:matrix_A}, \eqref{eq:matrix_B} and  \eqref{eq:vector_a}, respectively.
This problem is indeed in the form of \eqref{eq:mpec_prob}, and therefore, we can use Proposition \ref{pro:dc_decomp1} to obtain its DC formulation.

\subsubsection{Optimization over the solution-set of a variational inequality}
\label{subsubsec:OPoverVIP}
Let us consider a particular case of Problem \eqref{eq:mpec_prob} when the variable $y$ is absent. In this case, Problem \eqref{eq:mpec_prob}, with $S = \R^{m+n}$, takes the form
\begin{align}\label{eq:prob_P2}
&\min f(x) \tag{$\textrm{P}_2$}\\
&\textrm{s.t.}~~ x\in C, ~(Ax+a)^T(v-x) \geq 0, ~\forall v\in C,\notag
\end{align}
where, as before, $f$ is a real valued convex function on $C$ and $\emptyset \not=C \subseteq \R^n$ is a closed convex set.
Problems over the solution-set of a pseudomonotone variational inequality have been studied in \cite{Kalashikov1996} (notions of pseudomonotonicity and monotonicity can be taken from \cite{Konnov2000} or \cite{Konnov2001}). Here, we do not require any assumption on monotonicity. Note that without monotonicity of $A$, the solution-set of the variational inequality constraint in \eqref{eq:prob_P2} is not necessarily convex. Therefore, this problem remains a nonconvex optimization one. By Lemma \ref{le:gap_func} we can rewrite \eqref{eq:prob_P2} as
\begin{equation*}
\min\big\{f(x) ~|~ x\in C,~ g(x) \leq 0\big\}, 
\end{equation*}
where, by \eqref{eq:gap_func},
\begin{align*}
g(x)  = x^TAx -\frac{1}{2}x^TGx + a^Tx + \max_{v\in C}\big\{ -v^TAx - a^Tv - \frac{1}{2}v^TGv + x^TGv \big\}. 
\end{align*}
If $A$ is symmetric, we express $A$ as $A = A_1 - A_2$ with $A_1$ symmetric positive definite and $A_2$ symmetric positive (semi)definite. From Proposition \ref{pro:dc_decomp1} we have
\begin{align}\label{eq:g1_vip}
g_1(x) & = a^Tx + \max_{v\in C}\big\{ [(A_1+A_2)x - a]^Tv - v^TA_1v \big\},
\end{align}
and
\begin{equation}\label{eq:h1_vip}
h_1(x,y) = x^TA_2x.
\end{equation} 
Note that when $f $ is   constant, Problem \eqref{eq:prob_P2} becomes the affine variational inequality \cite{Facchinei2003, Lee2005}:
\begin{equation*}
\text{Find} ~ x \in C~~\text{such that:}~~ (Ax + a)^T(v-x) \geq 0, ~~\text{for all}~~ v\in C. 
\end{equation*}
By Lemma \ref{le:gap_func}, $x$ is a solution to this problem if and only if it is a global optimal solution to the differentiable DC program:
\begin{equation*}
0= \min\{ g(x):= g_1(x) - h_1(x): \ x\in C\}, 
\end{equation*}
where $g_1$ and $h_1$ are given as in Propositions \ref{pro:dc_decomp1}. 

\section{On Global Optimization Methods for AMPEC}
\label{sec:global_algorithm}
Theoretically, the global optimization methods such as branch-and-bound, outer and inner approximations, e.g., \cite{Horst1992}, can be applied to AMPEC by using the DC formulations obtained in the preceding sections. Note that AMPEC can be equivalently converted it into an one-level mathematical program with an additionally complementarity constraint by applying the Krush-Kuhn-Tucker theorem to the lower variational inequality. Branch-and-Bound algorithms have been developed in \cite{An2009, Anh2006, Muu2007, Thoai2005} for globally solving the latter problem. These existing algorithms use different subdivisions, but all of them take place in a space whose dimension is equal to the number of the Lagrangian multipliers. The latter number is large when the feasible set of the lower affine variational inequality is given, as usual, as $C :=\{ x\in R^n ~|~ x \geq 0,~ Px = q\}$ with $n$ large (often in practical problems). However, it is well recognized that global optimization algorithms work well only in the case when the dimension of the space, where the global optimization operations such as subdivision take place, is relatively small.

It can be observed that in AMPEC problem \eqref{eq:mpec_prob}, where $A$ is monotone on $C$, only the variable $y$ makes nonconvexity of the problem. In fact, when $A$ is monotone and $y$ is absent, the solution-set of the lower variational inequality is convex. This observation suggests us to look for DC decompositions of $g$ where the second component $h_1$ that makes $g$ nonconvex depends upon only $y$. From \eqref{eq:h1_func} in Proposition \ref{pro:dc_decomp1} and Remark \ref{re:B_expression} we see that if we choose $U = [(\Sigma B^{+})^T]^{+}$ and $A_2$ such that $2A_2 + U^TU = 0$, then $h_1$ is independent of $x$ and separable. In some models such as bilevel strongly convex quadratic problem \cite{Muu2003} and Nash-Cournot equilibrium model (Example 2.2.c), since $A$ is positive definite, one can choose $A_2 = - (1/2)U^TU$. Then, by virtue of Proposition \ref{pro:dc_decomp1}, we have $h_1(x,y) = \frac{1}{2} \norm{\Sigma y}^2$ is independent of $x$ and separable.

As an example, we now describe a branch-and-bound algorithm for minimizing a convex function over the equilibrium set of the Nash-Cournot equilibrium market model that we have studied in Subsection \ref{subsec:special_cases}. In practical Nash-Cournot models, the number $m$ of the materials that the producers need to produce the goods is much less than the number $n$ of the firms, for example, in electricity production, it takes only oil and coal as two main materials into account. 

This fact suggests that we should choose a DC decomposition such that the function $h_1$, which makes the problem nonconvex, depends upon only $y$ variable. For this purpose we choose the DC decomposition given in Proposition \ref{pro:dc_decomp1} with 
\begin{equation}\label{eq:matrix_A1}
A_1 = \begin{bmatrix}
2\beta &\beta &\beta &\cdots&\beta\\
\beta &2\beta &\beta &\cdots&\beta\\
\cdots&\cdots&\cdots&\cdots&\cdots\\
\beta &\beta &\beta &\cdots&2\beta
\end{bmatrix}
- \frac{1}{2}U^TU.
\end{equation}
and 
\begin{equation}\label{eq:matrix_A2}
A_2 = -\frac{1}{2}U^TU.
\end{equation}
Note that since $\lambda_{\min}(A) = \beta > 0$, where $\lambda_{\min}(A)$ is the smallest eigenvalue of $A$, the matrix $A$ is positive definite. 
If we choose $\Sigma$ such that $\lambda_{\max}( \Sigma ) < \beta$, where $ \lambda_{\max}( \Sigma )$ is the largest eigenvalue of $ \Sigma $, then $A_1$ is still positive definite.
By Proposition \ref{pro:dc_decomp1}, one has
\begin{align}\label{eq:g1_func_def3}
g_1(x,y) = \frac{1}{2}\norm{Ux + Vy}^2 + a^Tx + \max_{v\in C}\big\{ \!-v^TA_1v \!+\! [(A-U^TU)x \!-\! By - a]^Tv \!\big\},
\end{align}
\begin{equation}\label{eq:h1_func_def3}
h_1(x,y) = h_1(y):= \frac{1}{2}\norm{\Sigma y}^2 ~~~( \textrm{separable and depending on}~ y ~\textrm{only}).
\end{equation} 
Thus computing global optimal Nash equilibrium pairs to the bilevel Nash-Cournot equilibrium market model presented in Subsection \ref{subsec:special_cases} leads to the problem
\begin{equation}\label{eq:NC}
\alpha_{*} := \min \big\{f(x,y) ~|~ x\in X,~ y\in Y,~ g_1(x,y)-h_1(y)\leq 0 \big\}, \tag{NC} 
\end{equation}
where $g_1$ and $h_1$ are given by \eqref{eq:g1_func_def3} and \eqref{eq:h1_func_def3} respectively.

The separability property of $h_1$ suggests us to use the convex envelope
of $h_1$ on the box (rectangle) $Y$ to compute lower bounds in the branch-and-bound algorithm to be described below. Moreover, since $h_1$ depends upon only variable $y\in Y$, one can use an adaptive rectangular bisection that takes place in the $y$-space only.

Now we are going to describe in details these bounding and branching operations.

\subsection{Bounding by the convex envelope}
\label{subsec:bounding_convex_envelope}
We recall \cite{Falk1969, Horst1992} that a function $l(y)$ is said to be the convex envelope of a function $q(y)$ on a convex set $Y$ if $l$ is convex on $Y$, $l(y) \leq q(y)$ for every $y \in Y$ and if $p(y)$ is a convex function on $Y$ such that $p(y) \leq q(y)$ for every $y\in Y$ then $p(y) \leq l(y)$ for every $y\in Y$.
In general, computing the convex envelope of a function on an arbitrary convex set, even polyhedron, is expensive. Fortunately, in our case, since $h_1$ given in \eqref{eq:h1_func_def3} is separable, concave, and $Y$ is a box, its convex envelope is an affine function that can be given explicitly (see, e.g.\cite{Falk1969}). Namely, suppose that $h_1(y) = \sum_{j=1}^m \xi_j y_j^2$, $(\xi \geq 0$). Let $l^R$ denote the convex envelope of $-h_1$ on the box
\begin{equation*}
R:= \{ y=(y_1,\dots,y_m)^T ~|~ a_j \leq y_j \leq b_j, \ j=1,\dots,m\} \subseteq Y. 
\end{equation*}
Then $l^R(y) = \sum_{j=1}^m l_j^R(y_j)$, where $l_j^R$ is the convex envelope of the univariable function $-\xi_jy_j^2$ on the interval $[a_j, b_j]$ ($j=1,\dots,m)$. The latter in turn is the affine function joining $a_j$ and $b_j$.

Let $\alpha(R)$ and $\beta(R)$ denote the optimal value of Problem \eqref{eq:NC} restricted on $R$ and the optimal value of its relaxed problem,  respectively, that is
\begin{align}
&\alpha(R) := \min \big\{f(x,y)~|~ x\in X,~ y\in R,~ g_1(x,y)- h_1(y)\leq 0 \big\}, \label{eq:NC_R}\tag{$\textrm{NC}_R$}\\
&\beta(R) := \min \big\{f(x,y) ~|~ x\in X,~ y\in R,~ g_1(x,y)+ l^R(y)\leq 0 \big\}. \label{eq:RNC_R}\tag{$\textrm{RNC}_R$}
\end{align}
Since $l^R(y) \leq - h_1(y)$ for every $y\in R$, we have $\beta(R) \leq \alpha(R)$.

\subsection{An adaptive rectangular bisection}\label{subsec:bisection}
It is clear that if $\beta(R) = \alpha(R)$ then the minimum of $f$ over the set $x\in X, y\in R, g_1(x,y)- h_1(y)\leq 0$ has been found. Otherwise, if $\beta(R) < \alpha(R)$ then there must exist at least one index $j$ such that $l_j^R(y_j^{*}) <- \xi_jy^{*2}_j$, where $y_j^{*}$ denotes $j^{\textrm{th}}$ entry of an optimal solution to the relaxed problem defining $\beta(R)$. Let $j_R$ be an index such that
\begin{equation*}
\delta(R):= -\xi_{j_R}y^{*2}_{j_R} - l_j^R(y^{*}_{j_R}) = \max_{1\leq j\leq m}\big\{ -\xi_jy^{*2}_j - l_j^R(y^{*}_j)\big\}.
\end{equation*}
Note that at the ends of each edge of the box $R$, the value of the function $-\xi_j y^2_j$ and of its convex envelope coincide. Thus $\delta(R) \neq 0$ implies that $y^{*}_{j_R}$ is not an endpoint of the edge $j_R$ of $R$.

Using $j_R$ and $y^{*}_{j_R}$ we bisect $R$ into two subboxes $R^{+}$ and $R^{-}$ by setting
\begin{equation}\label{eq:Rplus_def}
R^{+} := \{ y=(y_1,\dots,y_m)^T \in R ~|~ y_{j_R} \geq y^{*}_{j_R}\},
\end{equation}
\begin{equation}\label{eq:Rminus_def} 
R^{-} := \{ y=(y_1,\dots,y_m)^T \in R ~|~ y_{j_R} \leq y^{*}_{j_R}\}.
\end{equation}
Clearly, both $R^{+}$ and $R^{-}$ are not empty. For this bisection we have the following lemma whose proof can be found, e.g., in \cite{Muu2000}:

\begin{lemma}\label{le:global_convergence} 
Let $\{R_k\}$ be an infinite sequence of boxes generated by the bisection process defined by \eqref{eq:Rplus_def} and \eqref{eq:Rminus_def}. Suppose that $R_{k+1} \subset R_k$ for every $k$. Then
\begin{equation*}
\lim_{k\to \infty} \big\{\alpha(R_k) - \beta(R_k)\big\} = 0. 
\end{equation*}
\end{lemma}

\subsection{Computing an upper bound}\label{subsec:upper_bound} 
Note that a feasible point of the AMPEC problem \eqref{eq:mpec_prob} can be computed whenever the lower problem is solved. In the Nash-Cournot equilibrium market model described in Subsection \ref{subsec:special_cases}, the lower problem can be solved efficiently with available codes, since it is a strongly convex quadratic program over the polyhedron $X$. 
In fact, with a fixed $y\in Y$, the lower problem is the strongly monotone variational inequality
\begin{equation}\label{eq:VIY_prob}
\text{Find}~ x\in X ~~\text{such that:}~~ (Ax +By +a)^T(v -x) \geq 0, ~~\text{for all}~ v\in X, \tag{$\textrm{VI}_y$} 
\end{equation}
where $A$ is given by \eqref{eq:matrix_A} and $a = (-\alpha,\dots,-\alpha)^T$.
This variational inequality is reduced to the strongly convex quadratic program (see, e.g. \cite{Konnov2000}):
\begin{equation*}
\min\big\{ \frac{1}{2}x^T A x + \sum_{k=1}^n(\mu_k+\alpha)x_k ~|~ x\in X \big\}, 
\end{equation*}
where $\mu_k = \sum_{i=1}^m c_{ki}y_i $. Hence, if $x$ is the optimal solution to this program then $(x,y)$ is a feasible point to the model, and therefore, $f(x,y)$ is an upper bound for the optimal value $\alpha_{*}$.
Now we are available to describe in detail an algorithm for global solving Problem \eqref{eq:NC} thereby obtaining a global optimal equilibrium pair to the bilevel Nash-Cournot equilibrium market model presented in Subsection \ref{subsec:special_cases}.

The B\&B algorithm is described as follows:

\noindent\textsc{B\&B Algorithm. }\\
\noindent{\bf Initialization. } Choose a tolerance $\varepsilon \geq 0$, take $R_0 = Y$ and solve the relaxed problem \eqref{eq:RNC_R} with $R = R_0$ to obtain the optimal value $\beta_0 :=\beta(R_0)$ and an optimal solution $(x^{R_0},y^{R_0})$. If $l_{R_0}(y^{R_0}) = h_1(y^{R_0})$, we are done: $(x^{R_0},y^{R_0})$ is a global optimal solution to Problem \eqref{eq:NC}. Otherwise, solve the lower problem \eqref{eq:VIY_prob} with $y= y^{R_0}$ to obtain a feasible point.  Let $(x^0,y^0)$ be the currently best feasible point and $\alpha_0 = f(x^0,y^0)$ be the currently best upper bound (we also call it the score). Set
\begin{align*}
\Gamma_0:= \begin{cases} 
\{R_0\} ~~&\text{if} ~ \alpha_0 - \beta_0 > \varepsilon (|\alpha_0|+1),\\
\emptyset ~~ &\text{otherwise}. 
\end{cases} 
\end{align*}
\noindent{\bf Iteration $k$} ($k=0,1,\dots)$. At the beginning of each iteration $k$ we have a family $\Gamma_k$ of subboxes of $Y$, a lower bound $\beta_k$, an upper bound $\alpha_k$ for the optimal value $\alpha_{*}$ and a feasible point $(x^k,y^k)$ such that $\alpha_k = f(x^k,y^k)$.
\begin{itemize}
\item[] a) If $\Gamma_k = \emptyset$, then terminate: $\alpha_k$ is an $\varepsilon$-solution and $(x^k,y^k)$ is an $\varepsilon$-global optimal solution.
\item[] b) If $\Gamma_k \neq \emptyset$, choose $\R_k \in \Gamma_k$ such that
\begin{equation*}
\beta(R_k) = \min\{ \beta(R)~|~ R \in \Gamma_k\}. 
\end{equation*}
Bisect $R_k$ into two rectangles $R_{k1}$ and $R_{k2}$ according to the bisection \eqref{eq:Rplus_def} and \eqref{eq:Rminus_def}. For each $(j=1, 2)$, compute
\begin{equation}\label{eq:RNC_R_kj}
\beta(R_{kj}):=\min \big\{f(x,y) ~|~ x\in X,~ y\in R_{kj}, ~ g_1(x,y) + l^{R_{kj}}(y)\leq 0 \big\} \tag{$\textrm{RNC}_{R_{kj}}$}. 
\end{equation}
\end{itemize}
Let $(x^{R_{kj}},y^{R_{kj}})$ be the obtained optimal solution to this subproblem. Use $y^{R_{kj}} $ $(j=1,2)$ to compute new feasible points by solving the strongly monotone variational inequalities \eqref{eq:VIY_prob} with $y= y^{R_{kj}}$ ($j=1,2)$.
 Let $(x^{k+1},y^{k+1})$ be the currently best feasible point and $\alpha^{k+1} = f(x^{k+1},y^{k+1})$ be the new upper bound (new score).
Delete all $ R\in \Gamma_k$ such that
\begin{equation*}
\alpha_{k+1} - \beta(R) \leq \varepsilon (|\alpha_{k+1}| + 1). 
\end{equation*}
Let $\Gamma_{k+1}$ the remaining set of subrectangles (may be empty). Then go to iteration $k$ with $k := k+1$. \hfill$\square$

Using Lemma \ref{le:global_convergence} by a standard argument commonly used in global optimization we can prove the following convergence property of the B \&B algorithm.

\begin{theorem}\label{th:global_convergence} Suppose that the sequence $\{(x^k, y^k)\}_k$ is generated by the B\& B  algorithm. Then
\begin{itemize}
\item[$\text{(i)}$] If the algorithm terminates at some iteration $k$ then $(x^k,y^k)$ is an $\varepsilon$- global optimal equilibrium pair to the Nash-Cournot equilibrium market model.
\item[$\text{(ii)}$] If the algorithm does not terminate then $\alpha_k \searrow \alpha_{*}$, $\beta_k \nearrow \alpha_{*}$ as $k \to +\infty$ and any limit point of the sequence $\{(x^k,y^k)\}$ is a global optimal equilibrium pair to the model.
\end{itemize}
\end{theorem}

\section{Numerical Results}\label{sec:numerical_results}
We have tested    the proposed  algorithm on the bilevel Nash-Cournot equilibrium market model presented in Subsection \ref{subsubsec:Nash_EP} with randomly generated data.
All computational results have been done in Matlab 7.8.0 (R2009a) for Linux running on a PC Desktop Intel(R) Core(TM)2 Quad CPU Q6600 2.4GHz, 3Gb RAM.
We generate data, choose parameters and solve the subproblems in the algorithm as follows:
\begin{itemize}
\item The objective function is chosen by a convex quadratic form $f(x,y) = \frac{1}{2}x^TQ_1x + \frac{1}{2}y^TQ_2y + q_1^Tx + q_2^Ty$, where $Q_1, Q_2, q_1$ and $q_2$ are generated randomly. The parameters $\beta=0.125$, $\alpha=10$ whereas $B = (c_{ij})_{n\times m}$ is generated randomly in $(0,1)$. The convex sets $X=[0, 5]^n$ and $Y = [0, 5]^m$, 
\item For computing the lower bound, we used the interior point method of the built-in Matlab solver FMINCON with maximum of iterations $500$ to solve the convex subproblems. The convex quadratic problems are solved by QUADPROG (a built-in Matlab solver) and CVX software (a freely available Malab code for convex programming). 
\item For computing the upper bound, a local optimization method in DC optimization is used that proves a feasible point to the problem \eqref{eq:mpec_c2}.
\end{itemize}
\begin{center}
\begin{tabular}{|c|c|c|c|c|c|r|c|l|}\hline\hline
\multicolumn{3}{|c|}{Problem Info.} & \multicolumn{6}{c|}{ Branch \& Bound algorithm}\\ \hline 
$\textrm{N}^0$ & $\textrm{m}$ & $\textrm{n}$ & \texttt{cbval} & \texttt{lbval} & \texttt{iter} & \texttt{time(s)} & \texttt{node} & \texttt{status}\\\hline
1 & 5 & 10 & 1338.2220 & 1338.2021 & 17 & 88.46 & 7 &  \texttt{solved} \\ \hline
2 & 10 & 10 & 1576.4746 & 1576.4407 & 154 & 962.99 & 56 &  \texttt{solved} \\ \hline 
3 & 5 & 20 & 3711.1289 & 3711.1289 & 43 & 521.23 & 7 &  \texttt{solved} \\ \hline
4 & 5 & 30 & 3537.2899 & 3537.2899 & 46 & 693.88 & 7 &  \texttt{solved} \\ \hline
5 & 8 & 50 & 3994.0027 & 3992.7705 & 99 & 7944.81 & 24 &  \texttt{solved} \\ \hline
6 & 5 & 100 & 3162.2176 & 3160.5017 & 47 & 7004.01 & 8 &  \texttt{solved} \\ \hline
7 & 6 & 100 & 4073.9795 & 4049.4880 & 62 & 9822.34 & 12 &  \texttt{incomp.} \\ \hline
8 & 7 & 100 & 3825.4430 & 3825.2157 & 73 & 11194.71 & 17 &  \texttt{solved} \\ \hline
9 & 5 & 150 & 2731.9005 & 2692.1867 & 43 & 14730.51 & 8 &  \texttt{incomp.} \\ \hline
10 & 6 & 150 & 3781.3484 & 3711.8269 & 73 & 20531.79 & 14 &  \texttt{incomp.} \\ \hline
11 & 1 & 200 & 3173.2954 & 3173.2954 & 9 & 4662.39 & 2 &  \texttt{solved} \\ \hline
12 & 2 & 200 & 2738.1198 & 2738.1198 & 19 & 5123.47 & 6 &  \texttt{solved} \\ \hline
13 & 3 & 200 & 2391.6111 & 2391.6111 & 18 & 6089.43 & 4 &   \texttt{solved} \\ \hline
14 & 4 & 200 & 2869.7684 & 2869.7684 & 22 & 10175.95 & 4 &  \texttt{solved} \\ \hline
15 & 5 & 200 & 3726.2399 & 3726.2399 & 55 & 26477.86 & 9 &  \texttt{solved} \\ \hline
16 & 6 & 200 & 2759.8484 & 2751.9396 & 75 & 36107.07 & 14 & \texttt{exceed} \\ \hline
17 & 7 & 200 & 2459.9965 & 2390.6909 & 78 & 36270.34 & 21 & \texttt{exceed} \\ \hline
18 & 8 & 200 & 3333.2645 & 3102.4295 & 80 & 36456.48 & 34 & \texttt{exceed} \\ \hline 
19 & 2 & 300 & 3008.2311 & 2975.1594 & 14 & 14963.37 & 2 & \texttt{incomp.} \\ \hline
20 & 3 & 300 & 3275.0818 & 3275.0818 & 29 & 29976.60 & 6 &  \texttt{solved} \\ \hline
\end{tabular}
\vskip 0.2cm
\centerline{Table 1. Computational results of B\&B algorithm for Nash-Cournot Problem}\label{tb:table1}
\end{center}

We perform the B\&B algorithm for $20$ random problem with different sizes. The results are reported in Table \ref{tb:table1}, where $m, n$ are the sizes of the problem; \texttt{iter} is the number of iterations; \texttt{cbval} is the currently best upper bound (score); \texttt{lbval} is the lower bound for the optimal value; \texttt{cputime} is the CPU time in second; \texttt{status} is the status of stopping criterion (\texttt{solved} shows that an $\varepsilon$- global optimal solution is found, \texttt{incomp.}  indicates that the program is stopped when the lower bound is improved too slowly, \texttt{exceed} means that the running time exceeds the limit $36.000$ seconds); and \texttt{node} is the maximum number of the nodes in the B\&B tree that have been stored. 

From the computational results we can observe the following technical remarks:
\begin{enumerate}
\item The proposed B\&B algorithm can solve globally AMPEC, in particular, bilevel convex quadratic problems, with several hundreds of decision variables while the number of the parameters is relatively small.
\item The numbers of iterations in Table \ref{tb:table1} indicates that the adaptive rectangular bisection used is effective.
\item Almost of the running time spends to solve the general convex subproblems for computing lower and upper bounds. Note that at each iteration in the interior point algorithms for convex subproblems one needs to solve strongly convex quadratic programs.
\end{enumerate}

\section{Conclusion}\label{sec:conclusion} 
We have formulated some classes of bilevel programming in forms of AMPEC. We have also used a regularization technique to obtain smoothly DC optimization formulations to AMPEC. 
A suitable regularization matrix results in a DC decomposition, where the second component depends upon only the parameter of the lower problem. We have described a decomposition branch-and-bound algorithm for globally solving AMPEC. This algorithm uses an adaptive rectangular bisection involving only the parameter which is often much less than the number of the decision variables in practical problems. Computational results on  a bilevel Nash-Cournot equilibrium market model show efficiency of the proposed algorithm.


\end{document}